\documentclass[12pt,a4paper]{amsart}

\usepackage{amscd}
\usepackage{graphicx}
\usepackage{amssymb,amsfonts,latexsym}
\usepackage{geometry}

\input xy
\xyoption {all}

\newtheorem{thm}{Theorem}[section]

\theoremstyle{definition}

\theoremstyle{remark}

\numberwithin{equation}{section}
\newcommand{\dQ}{\mathbb{Q}}
\newcommand{\dC}{\mathbb{C}}

\newcommand{\dZ}{\mathbb{Z}}
\newcommand{\dF}{\mathbb{F}}

\long\def\forget#1\forgotten{}

\begin{document}

\author{Joseph Cohen}
\address{
Department of Mathematics\\
Technion --- Israel Institute of Technology\\
Haifa, 32000\\
Israel }
\email{coheny@tx.technion.ac.il}


\author{Jack Sonn}
\address{
Department of Mathematics\\
Technion --- Israel Institute of Technology\\
Haifa, 32000\\
Israel }
\email{sonn@math.technion.ac.il}

\title[On $GCD(\Phi_N(a^n),\Phi_N(b^n))$]{On $GCD(\Phi_N(a^n),\Phi_N(b^n))$}
\date{\today}

\begin{abstract}
There has been interest during the last decade in properties of the sequence $\gcd(a^n-1,b^n-1), \ \
n=1,2,3,....$ where $a,b$ are fixed (multiplicatively independent) elements in one of $\dZ, \dC[T]$, or
$\dF_q[T]$ .  In the case of $\dZ$, Bugeaud, Corvaja and Zannier have obtained an upper bound $\exp (\epsilon n)$ for any given $\epsilon >0$ and all large $n$, and demonstrate its sharpness by extracting from a paper of Adleman,  Pomerance, and Rumely a lower bound $\exp(\exp(c\frac{\log
n}{\log\log n}))$  for infinitely many $n$, where $c$ is an absolute constant.  This paper generalizes these results to $\gcd(\Phi_N(a^n), \Phi_N(b^n))$ for any positive integer $N$, where $\Phi_N(x)$ is the $N$th cyclotomic polynomial, the preceding being the case $N=1$.  The upper bound follows easily from the original, but not the lower bound, which is the focus of this paper.
 The lower bound has been proved in the first author's Ph.D. thesis for the case $N=2$, i.e. for $\gcd(a^n+1,b^n+1)$.  In this paper we prove the lower bound for arbitrary $N$ under GRH (the generalized Riemann Hypothesis).  The analogue of the lower bound  for $\gcd(a^n-1,b^n-1)$ over $\dF_q[T]$ was proved by Silverman; we prove a corresponding generalization unconditionally.
\end{abstract}
\maketitle
\section{Introduction}\label{sec:intro}
In recent  years there has been interest \cite {BCZ},\cite {AR},\cite {si} \  in sequences of the form $$\gcd(a^n-1,b^n-1), \ \
n=1,2,3,....$$ where $a,b$ are fixed elements in one of $\dZ, \dC[T]$, or
$\dF_q[T]$.
Motivated by recurrence sequences and the Hadamard quotient theorem,  Bugeaud, Corvaja and Zannier \cite
{BCZ},   bounded the cancellation in the sequence $\frac{b^n-1}{a^n-1}$ by proving the following upper bound result:
\vskip 0.5em
\begin{thm}\label{thm:upper bound}\cite {BCZ} Let $a,b$ be multiplicatively independent
positive integers, $\epsilon>0$.  Then $\log\gcd(a^n-1,b^n-1)< \epsilon n$ for all sufficiently large $n$.
\end{thm}
Moreover, it is conjectured in \cite {AR} \ that if the  additional (necessary) condition $\gcd(a-1,b-1)=1$ holds, then $\gcd(a^n-1,b^n-1)=1$ for infinitely many $n$.

Returning to \cite {BCZ}, in order to show that Theorem  \ref{thm:upper bound} is close to best possible, it is remarked in \cite {BCZ} that one can derive  from a paper of Adleman,  Pomerance, and
Rumely \cite {APR} a lower bound result:
\begin{thm}\label{thm:lower bound}  \cite {BCZ} For any two positive integers $a,b$, there
exist infinitely many positive integers $n$ for which  $\log\gcd (a^n-1,b^n-1)>\exp(c\frac{\log
n}{\log\log n})$, where $c$ is an absolute constant.\rm \end{thm}
The result in \cite {APR} from which this is derived is an improvement of a result of Prachar \cite {pr}:
\begin{thm}\label{prachar} \cite{pr} Let $\delta(n)$ denote the number of divisors of $n$ of the form $p-1$, with $p$ prime.  Then there exist infinitely many $n$ such that $\delta(n)> \exp(c \log n/(\log \log n)^2).$ \end{thm}
The improvement in \cite {APR} (with a similar proof) removes the exponent 2 (and the $p-1$ are squarefree):
\begin{thm}\label{APR} \cite{APR} Let $\delta(n)$ denote the number of divisors of $n$ of the form $p-1$, with $p$ prime and $p-1$ squarefree.  Then there exist infinitely many $n$ such that $\delta(n)> \exp(c \log n/\log \log n).$ \end{thm}
It is interesting to note that in \cite{pr}, Prachar was motivated by a paper of N\"{o}bauer  \cite{no} which dealt with the group of invertible polynomial functions on $\dZ/n\dZ$ and particularly the subgroup of functions of the form $x^k$, whereas in \cite {APR}, Adleman, Pomerance and Rumely were motivated by the computation of a lower bound on the running time of a primality testing algorithm.
\vskip 0.5em
In his Ph.D. thesis \cite {coh}, \  the first author tests the robustness of these results and asks what happens to Theorem \ref{thm:lower bound} if $\gcd(a^n-1,b^n-1)$ is replaced by $\gcd(a^n+1,b^n+1)$ or by $\gcd(a^n+1,b^n-1)$, and proceeds to prove the analogous results for these sequences, using \cite {APR}:
\begin{thm}\label{thm:yossi}  \cite {coh} For any two positive nonsquare integers $a,b$, there
exist infinitely many positive integers $n$ for which  $\log\gcd (a^n+1,b^n+1)>\exp(c\frac{\log
n}{\log\log n})$, where $c$ is a constant depending on $a$ and $b$. The same result holds for $\gcd (a^n+1,b^n-1).$ \end{thm}

  (The corresponding analogues of Theorem  \ref{thm:upper bound} follow immediately from

  $x^n\pm1|x^{2n}-1$.)

If one  observes that the polynomials $x-1$ and $x+1$ are the first and second cyclotomic polynomials $\Phi_N(x)$, $N=1,2$, we ask if Theorems \ref{thm:upper bound} and \ref{thm:lower bound} also hold for $\gcd(\Phi_N(a^n), \Phi_N(b^n))$ for any positive integer $N$, or even for  $\gcd(\Phi_M(a^n), \Phi_N(b^n))$ for suitable positive integers $M,N$.  For Theorem \ref{thm:upper bound}, this is immediate
from $\Phi_N(x)|x^N-1$.  In this paper we deal with this generalization for Theorems \ref{thm:lower bound} (and \ref{thm:yossi}).
\vskip 0.5em
It should be remarked  that Corvaja and Zannier have made far-reaching generalizations of Theorem \ref{thm:upper bound} in \cite {CZ}, in other directions.  
\vskip 0.5em
In Section 2 we prove the above generalization for $\gcd(\Phi_N(a^n), \Phi_N(b^n))$ for any positive integer $N$ under the Generalized Riemann Hypothesis (GRH).  The explanation for this is that the generalization of Prachar's argument in this situation leads to an application of
the effective Chebotarev density theorem to a tower of Galois extensions $L_d/\dQ$,
where the exceptional zeros of the corresponding zeta functions of the $L_d$ are required to be bounded away from $1$ as $d$ goes to infinity \footnote{The 2-part of $[L_d:\dQ]$ is unbounded as $d\rightarrow \infty$ so the results of Stark and of Odlyzko and Skinner do not seem to apply.}.  Since we do not know if the exceptional zeros in our tower are bounded away from $1$, we  apply the stronger GRH version of the effective Chebotarev density theorem in which there are no exceptional zeros.  An additional attempt to avoid GRH using the Bombieri-Vinogradov theorem has so far not been successful. 
\vskip 0.5em
Silverman \cite {si} has proved an analogue of Theorem \ref{thm:lower bound} for the global function fields
$\dF_q(T)$:

\begin{thm}\label{thm:si} Let $\dF_q$ be a finite field and let
$a(T),b(T) \in \dF_q(T)$ be nonconstant monic polynomials.  Fix any
power $q^k$ of $q$ and any congruence class $n_0 +q^k\dZ \in
\dZ/q^k\dZ$.  Then there is a positive constant $c=c(a,b,q^k)>0$ such
that $$\deg (\gcd (a(T)^n-1,b(T)^n-1)) \geq cn$$ for infinitely many
$n \equiv n_0 \ (mod \ q^k)$. \end{thm}

In Section 3 we apply the method of Section 2 to prove (unconditionally) the corresponding cyclotomic generalization of Silverman's theorem.
\vskip 0.5em

\it Acknowledgment.  \rm We are grateful to Zeev Rudnick, Ram Murty and Jeff Lagarias for helpful discussions at various stages of the preparation of this paper. We also thank Joe Silverman for helpful comments on the initial draft.

\vskip 1em

\section{The case $a,b\in \dZ$ }\label{sec:rational case}
\vskip 1em
\begin{thm}[contingent on GRH]\label{thm:rational case}
Let $N$ be a
positive integer, $N=\ell_1^{s_1}\cdots \ell_r^{s_r},
\ell_1<\ell_2<\cdots \ell_r$, the factorization of $N$ into primes.
Let $a,b$ be positive integers, relatively prime to $N$, which are not $\ell_i$th powers in
$\dQ$ for $i=1,...,r$.
Then there exist infinitely many positive integers $n$ such that
$$\log\gcd (\Phi_N(a^n),\Phi_N(b^n))>\exp (\frac{c\log n}{\log \log n}),$$

where $c$ is a positive constant depending only on
$a,b,N$.
\end{thm}
\vskip 1em
\begin{proof} 
\ Suppose $p$ is a prime congruent to $1$ mod $N$ such that neither
$a$ nor $b$ is a $\ell_i$th power mod $p$ for $i=1,...,r$. Suppose
also that $n$ is a positive integer prime to $N$ and divisible by
$\frac{p-1}{N}$.  Then $p\mid{\text
{gcd}}(\Phi_N(a^n),\Phi_N(b^n))$.  Indeed, $(a^n)^N\equiv 1$ (mod
$p$).  The orders of $a$ and of $a^n$ mod $p$ are equal and divide
$N$.  If $a$ has order less than $N$, then there is a prime $\ell|N$
such that $a^{N/\ell}\equiv 1$ mod $p$, so $a^{(p-1)/\ell}\equiv 1$
mod $p$, whence $a$ is an $\ell$th power mod $p$, contrary to
hypothesis.  The idea of the proof of the theorem, a generalization of the proof in Prachar's paper, is to use the pigeonhole principle to produce, for large $x$, an $n \leq x^2$ with more than $\exp(c\frac{\log x}{\log \log x})$ divisors of the form $\frac{p-1}{N}$, $p$ prime, $c$ an absolute constant.  The result then follows.

Fix $0<\delta<1$.  Let $x$ be a positive real number and let
$K=K_{\delta}(x)$ be the product of all the primes $p\leq \delta \log x,$ $p
\nmid N$.   Let $A$ be the set of pairs $(m,p)$, $m$ a positive
integer, $p$ a prime,  $m \leq x$, $p \leq x$, $\gcd(m,N)=1$,  $p\equiv
1$ (mod $N$), $p \not\equiv 1$ (mod $N\ell_i$), $i=1,...,r$,  neither $a$
nor $b$ is an $\ell_i$th power mod $p$, $i=1,...,r$, and
$K|m\frac{p-1}{N}$.

Now for each $d|K$, let $A_d$ be the subset of $A$ consisting of
pairs $(m,p)\in A$ such that $(m,K)=K/d$
and $d|\frac{p-1}{N}$.  Let $N_0:=\ell_1\cdots \ell_r.$ We first
bound $|A_d|$ from below by bounding the following subset of $A_d$ of the form
$A_d'\times A_d''$, where
$$A_d'=\{m\leq x:(m,N_0K)=K/d\}$$ and
$$A_d''=\{p \leq x:p\equiv
1  \mod N,\ p \not\equiv 1  \mod N\ell_i,\ i=1,...,r, \ d\mid\frac{p-1}{N},$$
   $$ \text{ and neither } a \text{ nor } b \text{ is an } \ell_i\text{th power mod } p , \text{ } i=1,...,r \}.$$
 \vskip 0.5em
 To bound $|A_d' \times A_d''|$ from below, it suffices to bound each
of $|A_d'|,|A_d''|$ from below and take the product of the two lower
bounds.
 \vskip 0.5em
First, writing $d'=K/d$, $$|A_d'|=|\{m \leq
x:d'|m,(m/d',N_0K/d')=1\}|=|\{m/d' \leq x/d':(m/d',N_0K/d')=1\}|$$
$$\geq\phi(N_0K/d')[\frac{x/d'}{N_0K/d'}]=\phi(N_0d)[x/N_0K]$$ where $\phi$ denotes
Euler's $\phi$-function and $[-]$ the integer part.

To
bound $|A_d''|$ from below we use the effective form of Chebotarev's
density theorem due to Lagarias and Odlyzko \cite {LO} \ as
formulated by Serre \cite {Se} \  under the generalized Riemann
Hypothesis (GRH).

The condition $d\mid \frac{p-1}{N}$ is equivalent to $p\equiv 1$
(mod $Nd$), which is equivalent to $p$ splits completely in
$\dQ(\mu_{Nd})$, where $\mu_n$ denotes the group of $n$th roots of
unity.  The condition $a$ is an $\ell$th power mod $p$ ($\ell$
prime) is equivalent to  the condition $x^\ell-a$ has a root mod
$p$, which for $p\equiv 1$ modulo $\ell$ is equivalent to the
condition $x^\ell-a$ splits into linear factors mod $p$, which is
equivalent to the condition $p$ splits completely in (the Galois
extension) $\dQ(\mu_\ell,\root{\ell}\of{a})$ of $\dQ$, which for
$p\equiv 1$ (mod $Nd$) and $\ell\mid N$ is equivalent to $p$ splits
completely in $\dQ(\mu_{Nd},\root{\ell}\of{a})$.

Consider the Galois extension
$F_d=\dQ(\mu_{NdN_0},\root{N_0}\of{a},\root{N_0}\of{b)}$ of $\dQ$, with Galois group $G_d=G(F_d/\dQ)$, and the subset
$$C_d=G(F_d/\dQ(\mu_{Nd}))\setminus\{[\bigcup_i G(F_d/\dQ(\mu_{Nd},\root{\ell_i}\of{a}))]\bigcup [\bigcup_i G(F_d/\dQ(\mu_{Nd},\root{\ell_i}\of{b}))]\bigcup [\bigcup_i G(F_d/\dQ(\mu_{Nd\ell_i}))]\}$$ of $G_d$.  ($C_d$ is the complement in the first group you see in the display, of the union of all the other (sub)groups you see in the display.)

It is easily verified that $C_d$ is $G_d$-invariant under
conjugation, i.e. a union of conjugacy classes in $G_d$.
\vskip 0.5em
It follows from the definition of $C_d$ that

$$A_d''=\{p \leq x:p \text {  unramified in  } F_d, (p,F_d/\dQ)\subseteq C_d\}$$
where $(p,F_d/\dQ)$ denotes the Artin symbol. Set $$\pi_{C_d}(x):=
|A_d''|=|\{p \leq x: p \text {  unramified in  } F_d,
(p,F_d/\dQ)\subseteq C_d\}|.$$

By the effective Chebotarev density theorem cited above, under GRH
for the Dedekind zeta function of $F_d$,
$$R_d(x):=|\pi_{C_d}(x)-\frac{|C_d|}{|G_d|}Li(x)|\leq c_1\frac{|C_d|}{|G_d|}x^{1/2}(\log D_{F_d}+n_{F_d}\log x)$$
where $c_1$ is an absolute constant, $D_{F_d}$ is the discriminant
of $F_d$, $n_{F_d}=[F_d:\dQ]$ is the degree of $F_d$ over $\dQ$, and
$Li(x)$ is the logarithmic integral $\int_2^x \frac{dt}{\log t}$.
\vskip 0.5em
We have $|G_d|=\phi({Nd})\prod_i \ell_i^{e_i}$, where $e_i=2$ or $3$
according to whether or not $a,b$ are multiplicatively dependent mod
$\ell_i$th powers in $\dQ$.
\vskip 0.5em
\it Claim: $|C_d|=\prod_i (\ell_i-1)^{e_i}$.
\vskip 0.5em
Proof: \rm First we look at the case $r=1$ ($N$ is a power of $\ell_1$) and write $\ell=\ell_1$.  We need an elementary observation.
Let $H$ be the direct product of three cyclic groups of order $\ell$: $H=U\times V \times W$ with $U,V,W$ cyclic of order $\ell$.  Then $(u,v,w)\in
H\setminus ((U\times V)\cup (V\times W )\cup (U \times W)) \iff u\neq 1, v\neq 1, w\neq 1.$   Hence $|H\setminus ((U\times V)\cup (V\times W )\cup (U \times W))|=(\ell-1)^3$.
\vskip 0.5em
Now write $G(F_d/\dQ(\mu_{Nd}))\cong H_1
\times\cdots \times H_r$, where $H_i= U_i \times V_i\times W_i$ if $a,b$   are
multiplicatively independent mod $\ell_i$th powers in $\dQ$, and
$H_i= U_i \times  W_i$ if not.  The subgroups whose union
we are looking at (in the definition of $C_d$) can be identified with the subgroups of the form
$H_1 \times \cdots H_{i-1} \times U_i \times V_i \times H_{i+1} \times \cdots
\times H_r$ or $H_1 \times \cdots H_{i-1} \times V_i \times W_i \times H_{i+1}
\times \cdots \times H_r$ or $H_1 \times \cdots H_{i-1}\times U_i \times W_i  \times H_{i+1} \times \cdots
\times H_r$ for those $i$ for which $a,b$ are
multiplicatively \it independent \rm mod $\ell_i$th powers in $\dQ$,
and the subgroups $H_1 \times \cdots H_{i-1} \times U_i \times
H_{i+1} \times \cdots \times H_r$ or $H_1 \times \cdots H_{i-1} \times W_i \times
H_{i+1} \times \cdots \times H_r$ for those $i$ for which $a,b$ are
multiplicatively \it dependent \rm mod $\ell_i$th powers in $\dQ$.
An element $(h_1,...,h_r)$ is in the union of these $\iff$ some $h_i
\in (U_i\times V_i)\cup (V_i\times W_i )\cup (U_i \times W_i)$ for $i$ of the first kind, or $h_i\in U_i\cup W_i$ for some $i$
of the second kind. Hence $(h_1,...,h_r)$ lies in the complement (in $H$) of
the union $\iff h_i=(u_i,v_i,w_i)$ with $u_i,v_i,w_i\neq 1$ for all $i$ of the first
kind and $h_i=(u_i,w_i)$ with $u_i,w_i\neq 1$ for all $i$ of the second kind. It follows that
the complement has order $\prod_i ( \ell_i-1)^{e_i}$, proving the
claim.
\vskip 0.5em
We conclude that
$$\frac{|C_d|}{|G_d|}=\frac{\prod_i (\ell_i-1)^{e_i}}{\phi({Nd})\prod_i \ell_i^{e_i}}.$$

\vskip 0.5em
 By \cite {Se}, Prop. 5, p. 128,
$$\log D_{F_d} \leq (n_{F_d}-1) \sum_{p|Nabd}\log p + n_{F_d}\log n_{F_d}|\{p:p|Nabd\}|.$$
Now $$n_{F_d}=\phi({Nd})\prod_i
\ell_i^{e_i}\leq\phi(Nd)N^3=\phi(N)N^3\phi(d),$$ so
$$\log D_{F_d} \leq (\phi(N)N^3\phi(d)-1)\log^2(Nabd)+ \phi(N)N^3\phi(d)\log(\phi(N)N^3\phi(d))\log(Nabd)$$
$$\leq \phi(N)N^3\phi(d)\log^2(Nabd)+ \phi(N)N^3\phi(d)\log(\phi(N)N^3\phi(d))\log(Nabd)$$
$$ \leq 2(\phi(N)N^3ab)^3\phi(d)\log^2d$$
$$=f(N,a,b)\phi(d)\log^2d.$$

It now follows that
$$|A_d|\geq|A_d'\times A_d''|=|A_d'||A_d''|=|A_d'|\pi_{C_d}(x)$$
$$\geq \phi(N_0d)[\frac{x}{N_0K}](\frac{|C_d|}{|G_d|}Li(x)-c_1\frac{|C_d|}{|G_d|}x^{1/2}(\log D_{F_d}+n_{F_d}\log x))$$
$$\geq\phi(N_0d)[\frac{x}{N_0K}]\frac{|C_d|}{|G_d|}(Li(x)-c_1x^{1/2}(\log D_{F_d}+n_{F_d}\log x))$$
where $c_1$ is an absolute constant.  We now bound $$
Li(x)-c_1x^{1/2}(\log D_{F_d}+n_{F_d}\log x)$$ from below.  First,
$$\log D_{F_d}+n_{F_d}\log x \leq f(N,a,b)\phi(d)\log^2d + \phi(N)N^3\phi(d)\log x$$
$$\leq g(N,a,b)\phi(d) \log^2x\leq g(N,a,b)x^{\delta}\log x \leq g(N,a,b)x^{\delta+\epsilon}$$
(using $\phi(d)<d<K<x^{\delta}$ and $\log x < x^{\epsilon}$ for any
given $\epsilon$ and sufficiently large $x$). From this,
$$Li(x)-c_1x^{1/2}(\log D_{F_d}+n_{F_d}\log x) \geq \frac{x}{2\log x}-c_1x^{\frac{1}{2}+\delta +\epsilon}g(N,a,b)\geq
\frac{x}{4\log x}$$ (for sufficiently large $x$, using $Li(x) \sim
\frac{x}{\log x}$).  We then have
$$|A_d'|\pi_{C_d}(x)\geq \phi(N_0d)[\frac{x}{N_0K}]\frac{x}{4\log x}\frac{|C_d|}{|G_d|}$$
$$\geq \frac{1}{2}\phi(N_0d)\frac{x}{N_0K}\frac{x}{4\log x}\frac{\phi(N_0)^2}{\phi(N)N_0^2}\cdot\frac{1}{\phi(d)}$$
$$\geq \frac{1}{8}\frac{\phi(N_0)^3 x^2}{\phi(N)N_0^3 K \log x}=\frac{h(N)}{K}\frac{x^2}{\log x}.$$
It then follows that $$|A|=\sum_{d|K}|A_d|\geq
\frac{h(N)}{K}\frac{x^2}{\log
x}\sum_{d|K}1=\frac{h(N)}{K}\frac{x^2}{\log x}2^{\omega(K)}$$
$$\geq \frac{h(N)}{K}\frac{x^2}{\log x}2^{\frac{\frac{1}{4}\delta\log x}{\log \log x}}$$ where $\omega(K)$ denotes the number of primes dividing $K$.
For the last inequality we use  \cite {HW}, 22.2, p. 341, and 22.10,
p. 355:
$$ \omega(K)\sim \frac{\log K}{\log\log K}\Rightarrow \omega(K)\geq \frac{\log K}{2\log\log K}\geq \frac{1}{4}\frac{\delta\log x}{\log\log x}.$$
Now the number of positive integers $n\leq x^2$ such that $K|n$ is
at most $\frac{x^2}{K}$.  Furthermore, for every pair $(m,p)\in A$,
$m\frac{p-1}{N}$ is such an $n$.  Therefore there exists an $n\leq x^2$ such
that $K|n$ with at least
$$\frac{|A|}{x^2/K}>\frac{h(N)}{\log x}2^{\frac{\frac{1}{4}\delta\log x}{\log \log x}}=h(N)\exp(c_2\delta\frac{\log x}{\log\log x}-\log\log x)>\exp(c_3\frac{\log x}{\log\log x})$$
representations of the form $m\frac{p-1}{N}$, for $x$ sufficiently large,
where $c_2, c_3$ are absolute constants.  It follows that
$GCD(\Phi_N(a^n),\Phi_N(b^n))$ is a product of at least
$\exp(c_3\frac{\log x}{\log\log x})$ primes, hence is itself at
least $\exp\exp(c_4\frac{\log x}{\log\log x})$. As $n\leq x^2$ and
$\frac{\log x}{\log\log x}$ is an increasing function (for $x>e^e$),
the last expression is $\geq \exp\exp(c_5\frac{\log n}{\log\log
n})$.   \end{proof}
\vskip 1em

 The proof of Theorem 2.1 can be generalized to yield the following
\vskip 0.5em
\begin{thm} [contingent on GRH] \label{thm:genrationalcase}   Let $M,N$ be positive integers.  Let $D=\gcd(M,N)$ and assume $\gcd(M/D,D)=\gcd(N/D,D)=1$.  Let $L=lcm(M,N)=\ell_1^{s_1}\cdots \ell_r^{s_r},
\ell_1<\ell_2<\cdots \ell_r$, the factorization of $L$ into primes.
Let $a,b$ be positive integers, relatively prime to $L$, which are not $\ell_i$th powers in
$\dQ$ for $i=1,...,r$.
Then there exist infinitely many positive integers $n$ such that
$$\gcd (\Phi_M(a^n),\Phi_N(b^n))>(\exp (\exp (\frac{c\log n}{\log \log n}))),$$
where $c$ is a positive constant depending only on
$a,b,N$.\end{thm}
\vskip 0.5em
The proof is similar to the proof of Theorem \ref{thm:rational case}; we omit the details.  Also here, the case $M=1$, $N=2$ was proved unconditionally in \cite {coh}.
\vskip 0.5em

\vskip 2em

\section{The case  $a=a(T),b=b(T)\in\dF_q(T)$ }\label{sec:function field case}

\vskip 1em
In this section we will generalize Silverman's
 Theorem \ref{thm:si} \cite {si}:  Let $\dF_q$ be a finite field and let
$a(T),b(T) \in \dF_q(T)$ be nonconstant monic polynomials.  Fix any
power $q^k$ of $q$ and any congruence class $n_0 +q^k\dZ \in
\dZ/q^k\dZ$.  Then there is a positive constant $c=c(a,b,q^k)>0$ such
that $$\deg (\gcd (a(T)^n-1,b(T)^n-1)) \geq cn$$ for infinitely many
$n \equiv n_0 \ (mod \ q^k)$.

\vskip 1em
\rm The generalization will be as in the preceding section,
replacing $a(T)^n-1$ with $\Phi_m(a(T)^n)$ for an arbitrary fixed
positive integer $m$.  The proof will be similar in parts to the
proof of Theorem 1, but there will be some changes in notation.
\vskip 1em
\begin{thm} \label{function field case}  Let $\dF_q$ be a
finite field, and let $m$ be a positive integer prime to $q$,
$m=\ell_1^{e_1}\cdots \ell_s^{e_s}, \ell_1<\ell_2<\cdots \ell_s$,
the factorization of $m$ into primes. Let  $a(T),b(T) \in \dF_q(T)$
be nonconstant monic polynomials which are not $\ell_i$th powers in
$\dF_q(T)$ for $i=1,...,s$.  Fix  a power $q^k$ of $q$, and any
congruence class $n_0 +q^k\dZ \in \dZ/q^k\dZ$.  Then there is a
positive constant $c=c(m,a,b,q^k)>0$ such that $$\deg (\gcd
(\Phi_m(a(T)^n),\Phi_m(b(T)^n)) \geq cn$$ for infinitely many n
$\equiv n_0 \ (mod \ q^k)$.
\end{thm}

\vskip 1em
\it Proof. \rm Assume first that $(n_0,q)=1$.  Choose the smallest
positive integer $r$ such that $(r,m)=1$ and $rmn_0 \equiv -1 \ (mod
\ q^k).$    Let $Q=q^t$, where $t\geq k$ and  $q^t\equiv 1$ mod $mr$
(e.g. $t=k\phi(mr)$.   Let $n=\frac{Q^N-1}{mr}$, where $N$ is a
positive integer.  Let $\pi=\pi(T)$ be a monic irreducible
polynomial of degree $N$ in $\dF_Q[T]$ not dividing $a(T)b(T)$ (this
holds e.g. if $\deg(\pi)>\deg(a(T)b(T))$).  Then, writing
$a=a(T),b=b(T)$, $\pi|\Phi_m(a^n)$ if and only if $a^n$ is a
primitive $m$th root of unity mod $\pi$, i.e. $a^{nm}\equiv 1$ mod
$\pi$ and $a^{mn/\ell}\not\equiv 1$ mod $\pi$ for every $\ell|m$.
Substituting $n=\frac{Q^N-1}{mr}$, this holds $\Leftrightarrow$
$$a^{\frac{Q^N-1}{r}}\equiv 1 \ (mod \ \pi)$$ and
$$a^{\frac{Q^N-1}{r\ell}}\not\equiv 1 \ (mod \ \pi)$$ for all
$\ell|m$.  The first condition holds $\Leftrightarrow$ there exists
$A\in \dF_Q[T]$ such that $a\equiv A^r \ (mod \ \pi)$.  For such an
$A$, the second condition is equivalent to
$$A^{\frac{Q^N-1}{\ell}}\not\equiv 1 \ (mod \ \pi)$$ which is
equivalent to saying that $A$ is not an $\ell$th power mod $\pi$,
and since $(r,\ell)=1$, this is equivalent to saying that $a$ is not
an $\ell$th power mod $\pi$.  It follows that the two conditions
hold together $\Leftrightarrow$ $a$ is an $r$th power mod $\pi$ and
$a$ is not an $\ell$th power mod $\pi$ for all $\ell$ dividing $m$.
We conclude that $\pi|\Phi_m(a^n)$ if and only if $a$ is an $r$th
power mod $\pi$ and $a$ is not an $\ell$th power mod $\pi$ for all
$\ell$ dividing $m$.  Similarly, $\pi|\Phi_m(b^n)$ if and only if
$b$ is an $r$th power mod $\pi$ and $b$ is not an $\ell$th power mod
$\pi$ for all $\ell$ dividing $m$.
\vskip 0.5em
To count the number of $\pi$ dividing $\gcd
(\Phi_m(a^n),\Phi_m(b^n))$, we will use an effective version of
Chebotarev's density theorem for global function fields \cite {FJ},
p. 62, Prop. 5.16.  For this purpose, let
$$F:=\dF_{Q^N}(T)(\root{r}\of{a},\root{r}\of{b})$$ and let $$E:=\dF_{Q^N}(T)(\root{\ell_1}\of{a},\root{\ell_1}\of{b},...,\root{\ell_s}\of{a},\root{\ell_s}\of{b}).$$
\vskip 0.5em
Since $\deg\pi=N$, $\pi$ splits completely in $\dF_{Q^N}(T)$.
Therefore  $a$ and $b$ are $r$th powers mod $\pi$ if
and only if $\pi$ splits completely in $F$.  Furthermore, $a$ and $b$ are not $\ell$th powers mod $\pi$ for all
$\ell$ dividing $m$ if and only if $\pi$ does not split completely
in $\dF_{Q^N}(T)(\root{\ell}\of{a})$ nor in
$\dF_{Q^N}(T)(\root{\ell}\of{b})$ for all $\ell$ dividing $m$.
Accordingly, proceeding as in Section 2, consider the Galois
extension $EF/\dF_Q(T)$ with Galois group $G_N$, and let
$$C_N=G(EF/\dF_Q(T))\setminus \{[\bigcup_i G(EF/F(\root{\ell_i}\of{a}))]\bigcup [\bigcup_i G(EF/F(\root{\ell_i}\of{b}))]\}.$$
Then $\pi$ splits completely in $F$ and  $\pi$ does not split
completely in $\dF_{Q^N}(T)(\root{\ell}\of{a})$ nor in
$\dF_{Q^N}(T)(\root{\ell}\of{b})$ for all $\ell$ dividing $m$, if
and only if $(\pi, EF/\dF_Q(T))\subseteq C_N$.
\vskip 0.5em
Now the same counting argument as in the preceding section gives
$|G_N|=Nr^2\prod_i\ell_i^{e_i}$ and $|C_N|=\prod_i(\ell_i-1)^{e_i}$.
Applying \cite {FJ}, p. 62, Prop. 5.16 \footnote{This is an effective Chebotarev density theorem for global function fields, implied by the Riemann Hypothesis for curves over finite fields, which is a theorem.} (and observing that a
conjugacy class can be replaced by any union of conjugacy classes in
that theorem), we get
$$|\{\pi\in \dF_Q[T]:\pi \ \text{monic irreducible of degree} \  N, \ (\pi, EF/\dF_Q(T))\subseteq C_N\}|$$ $$=\frac{|C_N|}{|G_N|}Q^N+O(Q^{N/2})=\frac{\prod_i(\ell_i-1)^{e_i}}{Nr^2\prod_i\ell_i^{e_i}}Q^N+O(Q^{N/2}).$$  It follows that $$\deg (\gcd (\Phi_m(a(T)^n),\Phi_m(b(T)^n)) \geq N(\frac{\prod_i(\ell_i-1)^{e_i}}{Nr^2\prod_i\ell_i^{e_i}}Q^N+O(Q^{N/2}))$$ $$=
\frac{\prod_i(\ell_i-1)^{e_i}}{r^2\prod_i\ell_i^{e_i}}Q^N+O(NQ^{N/2})\geq
cn$$ for some constant not depending on $N$, and $n=Q^N-1$.  This
proves Theorem 2 when $(n_0,q)=1$.  The  case $(n_0,q)\neq 1$
follows from the case $(n_0,q)=1$ as in \cite {si}. \qed
\vskip 0.5em
As in the previous section, the proof of Theorem 3.1 can be generalized to yield the following

\begin{thm} \label{thm:genfunctionfieldcase} Let $\dF_q$ be a finite field, $u,v$ be positive integers, $d=\gcd(u,v)$, and assume $\gcd(u/d,d)=\gcd(v/d,d)=1$.  Let $a=a(T)$, resp.
$b=b(T)\in \dF_q[T]$ be monic nonconstant polynomials which are not $\ell$th powers in $\dF_q[T]$ for all $\ell|u$, resp. $\ell|v$.  Fix  a power $q^k$ of $q$, and any
congruence class $n_0 +q^k\dZ \in \dZ/q^k\dZ$.  Then there is a
positive constant $c=c(u,v,a,b,q^k)>0$ such that $$\deg (\gcd
(\Phi_u(a(T)^n),\Phi_v(b(T)^n)) \geq cn$$ for infinitely many $n
\equiv n_0 \ (mod \ q^k)$.\end{thm}
\vskip 0.5em
The details are omitted.

\vskip 2em
\bibliographystyle{plain}

\begin{thebibliography}{10}




 \bibitem{APR}
 L. M. Adleman,  C. Pomerance and R.S. Rumely,
 \newblock  On distinguishing prime numbers from composite numbers,
\newblock {\em Ann. Math.} 117:  173--206, 1983

\bibitem{AR}
N. Ailon and Z. Rudnick,
\newblock Torsion points on
curves and common divisors of $a^k-1,b^k-1$,
\newblock  {\em Acta Arith. }
113:31--38, 2004

\bibitem{BCZ}
Y. Bugeaud, P. Corvaja and U. Zannier,
\newblock An
upper bound for the G.C.D. of $a^n - 1$ and $b^n -1$,
\newblock  {\em Math.
Zeit. }
243:79--84, 2003

\bibitem{coh}
 J. Cohen,
 \newblock  Primitive roots in algebraic number fields,
\newblock {\em Ph.D. Thesis, Technion},  2004


\bibitem{CZ}
P. Corvaja and U. Zannier,
\newblock A lower bound for the height of a rational function at S-unit points,
\newblock {\em Monatsh. Math.} 144;3:203–-224, 2005

\bibitem{FJ}
M. Fried and M. Jarden,
\newblock {\em Field Arithmetic},
\newblock Springer-Verlag , New York-Heidelberg, 1986

\bibitem{HW}
G.H. Hardy and E.M. Wright,
\newblock {\em An Introduction
to the Theory of Numbers (Fifth Ed.)},
\newblock Oxford Univ. Press, Oxford, 1979

\bibitem{LO}
J. Lagarias and A.M. Odlyzko,
\newblock Effective
versions of the Chebotarev density theorem,
\newblock {\em Algebraic
Number Fields: L-functions and Galois properties (Proc. Sympos.,
Univ. Durham, 1975)}, Academic Press, London, 409-464, 1977

\bibitem{no}
W. N\"{o}bauer,
\newblock \"{U}ber eine Gruppe der Zahlentheorie,
\newblock {\em Monatsh. Math.}, 58:181-192 (1954)

\bibitem{pr}
K. Prachar,
\newblock \"{U}ber die Anzahl der Teiler einer nat\"{u}rlichen Zahl, welche die Form $p-1$  haben,
\newblock {\em Monatsh. Math.}, 59:91–-97 (1955)

\bibitem{Se}
J.-P. Serre,
\newblock Quelques applications du
theoreme de densite de Chebotarev,
\newblock {\em Publ. Math. IHES}, 54:123-201, 1982

\bibitem{si}
J. Silverman,
\newblock Common divisors of $a^n-1$ and
$b^n-1$ over function fields,
\newblock {\em New York J. Math.}, 10:37--43 , 2004

\end{thebibliography}

\def\cprime{$'$}

\end{document}